\documentclass[english,12pt]{smfart}
\usepackage{amssymb}
\usepackage{amsfonts}
\usepackage{amscd}
\usepackage[all]{xypic}
\usepackage[T1]{fontenc}

\swapnumbers

\def\deg{{\rm deg}}

\def\Var{{\rm Var}}
\def\hsp{{\rm hsp}}

\def\GL{{\rm GL}}
\def\SL{{\rm SL}}
\def\HS{{\rm HS}}

\let\cal\mathcal

\def\AA{{\mathbf A}}

\def\CC{{\mathbf C}}

\def\GG{{\mathbf G}}

\def\LL{{\mathbf L}}

\def\NN{{\mathbf N}}

\def\QQ{{\mathbf Q}}

\def\ZZ{{\mathbf Z}}

\def\cD{{\mathcal D}}

\def\cL{{\mathcal L}}
\def\cM{{\mathcal M}}

\def\cS{{\mathcal S}}

\def\cV{{\mathcal V}}
\def\cW{{\mathcal W}}
\def\cX{{\mathcal X}}
\def\cY{{\mathcal Y}}
\def\cZ{{\mathcal Z}}
\mathchardef\alphag="7C0B
\mathchardef\betag="7C0C
\mathchardef\gammag="7C0D
\mathchardef\deltag="7C0E
\mathchardef\varepsilong="7C22
\mathchardef\varphig="7C27
\mathchardef\psig="7C20
\mathchardef\zetag="7C10
\mathchardef\epsilong="7C0F
\mathchardef\rhog="7C1A
\mathchardef\taug="7C1C
\mathchardef\upsilong="7C1D
\mathchardef\iotag="7C13
\mathchardef\thetag="7C12
\mathchardef\pig="7C19
\mathchardef\sigmag="7C1B
\mathchardef\etag="7C11
\mathchardef\omegag="7C21
\mathchardef\kappag="7C14
\mathchardef\lambdag="7C15
\mathchardef\mug="7C16
\mathchardef\xig="7C18
\mathchardef\chig="7C1F
\mathchardef\nug="7C17
\mathchardef\varthetag="7C23
\mathchardef\varpig="7C24
\mathchardef\varrhog="7C25
\mathchardef\varsigmag="7C26
\mathchardef\Omegag="7C0A
\mathchardef\Thetag="7C02
\mathchardef\Sigmag="7C06
\mathchardef\Deltag="7C01
\mathchardef\Phig="7C08
\mathchardef\Gammag="7C00
\mathchardef\Psig="7C09
\mathchardef\Lambdag="7C03
\mathchardef\Xig="7C04
\mathchardef\Pig="7C05
\mathchardef\Upsilong="7C07

\newtheorem{theorem}[subsection]{Theorem}
\newtheorem{lem}[subsection]{Lemma}
\newtheorem{cor}[subsection]{Corollary}
\newtheorem{prop}[subsection]{Proposition}

\theoremstyle{definition}
\newtheorem{definition}[subsubsection]{Definition}

\newtheorem{def-prop}[subsubsection]{Definition-Proposition}

\theoremstyle{remark}
\newtheorem{remark}[subsection]{Remark}

\theoremstyle{plain}

\numberwithin{equation}{subsection}

\def\boxit#1#2{\setbox1=\hbox{\kern#1{#2}\kern#1}%
\dimen1=\ht1 \advance\dimen1 by #1
\dimen2=\dp1 \advance\dimen2 by #1
\setbox1=\hbox{\vrule height\dimen1 depth\dimen2\box1\vrule}%
\setbox1=\vbox{\hrule\box1\hrule}%
\advance\dimen1 by .4pt \ht1=\dimen1
\advance\dimen2 by .4pt \dp1=\dimen2 \box1\relax}

\renewcommand{\theequation}{\thesubsection.\arabic{equation}}

\let\cal\mathcal

\def\AA{{\mathbf A}}

\def\CC{{\mathbf C}}

\def\GG{{\mathbf G}}

\def\LL{{\mathbf L}}

\def\NN{{\mathbf N}}

\def\QQ{{\mathbf Q}}

\def\ZZ{{\mathbf Z}}

\def\cD{{\mathcal D}}

\def\cL{{\mathcal L}}
\def\cM{{\mathcal M}}

\def\cS{{\mathcal S}}

\def\cV{{\mathcal V}}
\def\cW{{\mathcal W}}
\def\cX{{\mathcal X}}
\def\cY{{\mathcal Y}}
\def\cZ{{\mathcal Z}}
\mathchardef\alphag="7C0B
\mathchardef\betag="7C0C
\mathchardef\gammag="7C0D
\mathchardef\deltag="7C0E
\mathchardef\varepsilong="7C22
\mathchardef\varphig="7C27
\mathchardef\psig="7C20
\mathchardef\zetag="7C10
\mathchardef\epsilong="7C0F
\mathchardef\rhog="7C1A
\mathchardef\taug="7C1C
\mathchardef\upsilong="7C1D
\mathchardef\iotag="7C13
\mathchardef\thetag="7C12
\mathchardef\pig="7C19
\mathchardef\sigmag="7C1B
\mathchardef\etag="7C11
\mathchardef\omegag="7C21
\mathchardef\kappag="7C14
\mathchardef\lambdag="7C15
\mathchardef\mug="7C16
\mathchardef\xig="7C18
\mathchardef\chig="7C1F
\mathchardef\nug="7C17
\mathchardef\varthetag="7C23
\mathchardef\varpig="7C24
\mathchardef\varrhog="7C25
\mathchardef\varsigmag="7C26
\mathchardef\Omegag="7C0A
\mathchardef\Thetag="7C02
\mathchardef\Sigmag="7C06
\mathchardef\Deltag="7C01
\mathchardef\Phig="7C08
\mathchardef\Gammag="7C00
\mathchardef\Psig="7C09
\mathchardef\Lambdag="7C03
\mathchardef\Xig="7C04
\mathchardef\Pig="7C05
\mathchardef\Upsilong="7C07

\def\ord{{\rm ord}}

\begin{document}

\title[Motivic zeta functions and
castling transformations]{Motivic
zeta functions for prehomogeneous vector spaces and castling transformations}

\author{Fran\c cois Loeser}

\address{{\'E}cole Normale Sup{\'e}rieure,
D{\'e}partement de math{\'e}matiques et applications,
45 rue d'Ulm,
75230 Paris Cedex 05, France
(UMR 8553 du CNRS)}
\email{Francois.Loeser@ens.fr}
\urladdr{http://www.dma.ens.fr/~loeser/}



\maketitle

\section{Introduction}
\subsection{}Let us recall the classification of irreducible regular prehomogeneous
vector spaces due to Sato and Kimura \cite{sato}.
An irreducible
regular prehomogeneous vector spaces over a field $k$ 
\cite{sato} 
consists of the datum $(G, X)$ of
a connected reductive $k$-algebraic group
$G$ 
together with a linear representation $\rho$ of $G$
on a finite dimensional  affine space
$X= \AA^m_k$ over $k$, such the action is transitive on the complement
of a geometrically irreducible hypersurface $f = 0$ in
$X$, called the singular locus.
In the theory of Sato and Kimura,
a fundamental role is played by castling transforms,
which are defined as follows.

Let $G$ be 
a connected reductive $k$-algebraic group 
and $\rho$ be a linear representation of $G$
on 
$\AA^m_k$.
Write
$m = r_1 + r_2$.
Then set $X_1 = M_{m, r_1}$ 
and let 
$G \times \SL_{r_1}$
act on $X_1$
by 
$$\rho_1 (g, g_1) x_1 = \rho (g) x_1 {}^{t}g_1.$$
Similarly
set $X_2 = M_{m, r_2}$
and let 
$G \times \SL_{r_2}$
act on $X_2$
by 
$$\rho_2 (g, g_2) x_2 = {}^t \rho (g)^{-1} x_2 {}^{t}g_2.$$
Assume $G_j = {\rm Im} \rho_j$
is an irreducible
regular prehomogeneous vector space for $j = 1, 2$ with
singular locus the irreducible hypersurface
$f_j = 0$.
In this case one says that
the prehomogeneous vector spaces $(G_1, X_1)$
and
$(G_2, X_2)$ are related by a castling transformation.
Furthermore, cf. \ref{tsc}, $f_1$ and $f_2$ are homogeneous
polynomials whose degrees satisfy
${\deg} f_j = r_j d$, for some integer $d$.
We shall say that two 
irreducible regular prehomogeneous
vector spaces
over $k$ are equivalent if, up to isomorphism,
they may be included in a finite chain of 
irreducible regular prehomogeneous
vector spaces whose consecutive terms are related by a castling
transform.
In each equivalence class there is a unique 
prehomogeneous
vector space with $\dim X$ minimal. Such an object is called reduced in
the terminology of 
Sato and Kimura. When $k$ is algebraically closed of characteristic zero,
the set of all reduced 
irreducible regular prehomogeneous
vector spaces has been divided into 29 types by Sato and Kimura
\cite{class}. If $k$ is of characteristic zero and one restricts to
$k$-split groups, then the absolute classification restricts to a relative one.

\subsection{}Thus, if one wants to compute the value of some invariant
for irreducible
regular prehomogeneous vector spaces, it is enough: a) to understand the behaviour of the invariant under castling transformations ;
b) to compute the invariant explicitely for reduced 
prehomogeneous vector spaces.  
For instance, when $k = \CC$,
if one denotes by $b_j$ the $b$-function of 
$f_j$, we have the following relation due to Shintani (cf. \cite{kimura}):
\begin{equation}
b_2 (s) \prod_{1 \leq i \leq r_1}
\prod_{0 \leq j < d} \Bigl(s + \frac{(i + j)}{d}\Bigr)
=
b_1 (s) \prod_{1 \leq i \leq r_2}
\prod_{0 \leq j < d}\Bigl(s + \frac{(i + j)}{d}\Bigr),
\end{equation}
and the $b$-functions of all 29 types,
except for one, have been tabulated by Kimura in \cite{kimura}.
When $k$ is a $p$-adic field and $Z_j$ denotes
Igusa's local zeta function attached to $f_j$, then Igusa proved in
\cite{igusa} the relation
\begin{equation}\label{ig}
Z_2 (s) \prod_{1 \leq j \leq r_1}
\Bigl(\frac{1 - q^{-j}}{1 - q^{- (ds + j)}}\Bigr)
=
Z_1 (s) \prod_{1 \leq j \leq r_2}
\Bigl(\frac{1 - q^{-j}}{1 - q^{- (ds + j)}}\Bigr),
\end{equation}
where $q$ denote the cardinality of the residue field
and the coefficients of each  polynomial $f_j$ are assumed
to belong to the valuation ring but not all to the valuation ideal ;
furthermore, Igusa was able to compute $Z (s)$
for 24 out of the 29 types, see 
\cite{aa}, \cite{igusaclass}.

\subsection{}
Another important invariant of singularities of functions
is the Hodge spectrum
\cite{St1}, \cite{St2}, \cite{Va} and a natural question
is to understand its behaviour under castling 
transformations. It is not clear to us how to answer this question
using classical techniques, and the aim of the present paper
is to explain how this can be achieved using motivic integration.
In fact, replacing $p$-adic integration with integration on arc spaces,
Denef and Loeser defined in
\cite{motivic}, \cite{lef}, \cite{barc} a motivic zeta function
$Z_f$ 
for a function $f$
which is a geometric analogue of Igusa's local zeta function.
The motivic zeta function
$Z_f$ contains a great amount of information on $f$. In particular one
can read off the whole Hodge spectrum of $f$ from
$Z_f$. 
The main result of the present paper is Theorem \ref{MT}
where a
relation completely analogous to (\ref{ig}) is proved for 
motivic zeta functions. We then deduce the
behaviour of the Hodge spectrum  under castling 
transformations
in Corollary \ref{sp}.
To prove Theorem \ref{MT}, it would be natural to
imitate Igusa's proof in the
$p$-adic case as much as possible.
In fact, we were not able 
to follow  Igusa's proof to closely, since in that proof
crucial use is made
of properties of the Haar measure on a locally compact
$p$-adic group, especially of its uniqueness, for which no analogue is
available yet on the motivic integration side, though highly desirable.
Hence we are led to follow a more down to earth approach,
relying on some direct computations in spaces of matrices, as in 
Lemma \ref{ML}.

Finally, in section \ref{gen} we generalize our results
to 
prehomogeneous vector spaces which are no longer assumed to be
irreducible and regular.

\section*{}I would like to thank Jan Denef and Akihiko Gyoja
for conversations that raised
my interest 
on the questions considered in the present paper.

\section{Notations and conventions}
\subsection{Grothendieck groups of varieties}We fix a field $k$. 
By a variety over $k$ we
shall mean a separated and reduced
scheme of finite type over $k$. 
Let $S$ be an algebraic variety over $k$.  By an $S$-variety we mean a
variety $X$ together with a morphism $X \rightarrow S$.  The $S$-varieties 
form
a category denoted by $\Var_S$, the arrows being
the morphisms that commute with
the morphisms to $S$.

We denote by $K_0 (\Var_S)$ the Grothendieck group of
$S$-varieties. It is
an abelian group generated by symbols $[X]$, for $X$ an $S$-variety, with the
relations $[X] = [Y]$ if $X$ and $Y$ are isomorphic in
$\Var_S$, and $[X] = [Y]
+ [X \setminus Y]$ if $Y$ is Zariski closed in $X$.  There is a natural ring
structure on $K_0 (\Var_S)$, the product of $[X]$ and $[Y]$ being equal to $[X
\times_S Y]$.  Sometimes we will also write $[X/S]$ instead of $[X]$, to
emphasize the role of $S$.  We write $\LL$ to denote the class of $\AA^1_k \times S$ in $K_0(\Var_S)$, where the morphism from $\AA^1_k
\times S$ to $S$ is the natural projection.  We denote by $\mathcal{M}_S$ the
ring
obtained from $K_0(\Var_S)$ by inverting $\LL$.  When $A$ is a
constructible subset of some $S$-variety, we define $[A/S]$ in the obvious way,
writing $A$ as a disjoint union of a finite number of
locally closed subvarieties $A_i$.  Indeed
$[A/S] := \sum_i [A_i/S]$ does not depend on the choice of the
subvarieties
$A_i$. 

When $S$ consists of only one geometric point, \textit{i.e.} $S = {\rm Spec} (k)$, then
we will write $K_0(\Var_k)$ instead of $K_0(\Var_S)$
and $\mathcal{M}_k$
instead of
$\mathcal{M}_S$. 

\subsection{Equivariant Grothendieck groups}
We need some technical preparation in order to take care of the monodromy
actions in the next section.

For any positive integer $n$, let $\mu_n$ be the group of all $n$-th roots of
unity (in some fixed algebraic closure of $k$).  Note that $\mu_n$'s
is actually
an algebraic variety over $k$, namely ${\rm Spec} (k[x]/(x^n-1))$.
The $\mu_n$ form a projective system, with respect to the maps $\mu_{nd}
\rightarrow \mu_n : x \mapsto x^d$.  We denote by $\hat \mu$ the projective
limit of the $\mu_n$.  Note that the group $\hat \mu$ is not an algebraic
variety, but a pro-variety.

Let $X$ be an $S$-variety.  A good $\mu_n$-action on $X$ is a group
action
$\mu_n \times X \rightarrow X$ which is a morphism of $S$-varieties, such that
each orbit is contained in an affine subvariety of $X$.  This last condition is
automatically satisfied when $X$ is a quasi projective variety.  A good
$\hat
\mu$-action on $X$ is an action of $\hat \mu$ on $X$ which factors through a
good $\mu_n$-action, for some $n$.

The monodromic Grothendieck group $K^{\hat \mu}_0
(\Var_S)$ is defined as
the
abelian group generated by symbols $[X, \hat \mu]$ (also denoted by $[X/S,\hat
\mu]$, or simply $[X]$), for $X$ an $S$-variety with good $\hat \mu$-action,
with the relations $[X, \hat \mu]  = [Y,\hat \mu]$ if $X$ and $Y$ are
isomorphic as $S$-varieties with $\hat \mu$ -action, and $[X,\hat \mu] = [Y,\hat
\mu] + [X \setminus Y, \hat \mu]$ if $Y$ is Zariski closed in $X$ with the
$\hat \mu$-action on $Y$ induced by the one on $X$, and moreover $[X \times
V,\hat \mu] = [X \times \AA^n_k, \hat \mu]$ where $V$ is the $n$-dimensional
affine space over $k$ with any linear $\hat \mu$-action, and $\AA^n_k$ is
taken with the trivial $\hat \mu$-action.  There is a natural ring structure on
$K^{\hat \mu}_0 (\Var_S)$, the product being induced by the fiber product over
$S$.  We write $\LL$ to denote the class in $K_0^{\hat \mu}
(\Var_S)$ of
$\AA^1_k \times S$ with the trivial $\hat \mu$-action.

We denote by $\mathcal{M}^{\hat \mu}_S$ the ring obtained from $K_0^{\hat
\mu}(\Var_S)$ by inverting $\LL$.  When $A$ is a constructible subset
of $X$ which is stable under the $\hat \mu$-action, then we define $[A,\hat
\mu]$ in the obvious way.  When $S$ consists of only one geometric point,
\textit{i.e.}
$S = {\rm Spec} (k)$, then we will write $K^{\hat \mu}_0(\Var_k)$
instead of $K^{\hat \mu}_0(\Var_S)$. 

Note that for any $s \in S(k)$ we have natural morphisms
${\rm Fiber}_s : K^{\hat \mu}_0 (\Var_S) \rightarrow K^{\hat \mu}_0 (\Var_k)$
and ${\rm Fiber}_s : \mathcal{M}^{\hat \mu}_S \rightarrow \mathcal{M}^{\hat \mu}_k$
given by
$[X,\hat
\mu] \rightarrow [X_s,\hat \mu]$, where $X_s$ denotes the fiber at $s$ of $X
\rightarrow S$.

We shall also consider the canonical morphism
$\mathcal{M}^{\hat \mu}_S \rightarrow \mathcal{M}^{\hat \mu}_k$
which to $[X/S,\hat \mu]$ assigns
$[X,\hat \mu]$.

\subsection{The arc space of $X$}
For each natural number $n$ we consider the space ${\cal L}_n(X)$ of arcs
modulo
$t^{n+1}$ on $X$. This is an algebraic variety over $k$, whose $K$-rational
points, for any field $K$ containing $k$, are the $K[t]/t^{n+1}K[t]$-rational
points of $X$.  For example when $X$ is an affine variety with equations $f_i
(x) = 0$, $i= 1, \cdots, m$, $x = (x_1,\cdots, x_n)$, then ${\cal L}_n(X)$
is given by the equations, in the variables $a_0, \cdots, a_n$,
expressing that $f_i (a_0 + a_1t + \cdots + a_n t^n) \equiv 0
\mod t^{n+1}, i = 1,\cdots, m$.  

Taking the projective limit of these algebraic varieties ${\cal L}_n(X)$ we
obtain
the arc space ${\cal L}(X)$ of
$X$, which is a reduced separated scheme
over
$k$.  In general, ${\cal L}(X)$ is not of finite type over $k$ 
(\textit{i.e.} ${\cal L}(X)$ is an ``algebraic
variety of infinite dimension'').  The $K$-rational points of ${\cal L} (X)$ are
the
$K[[t]]$-rational points of $X$.  These are called $K$-arcs on $X$.  For
example when $X$ is an affine complex variety with equations $f_i(x) = 0,
i = 1, \cdots m, x = (x_1,\cdots,x_n)$, then the $\CC$-rational points
of
${\cal L}(X)$ are the sequences $(a_0, a_1, a_2, \cdots) \in
(\CC^n)^{\NN}$ satisfying $f_i (a_0 + a_1 t + a_2t^2 +
\cdots )  = 0$,  for $i = 1,\cdots,m$. 
For any $n$, and for $m > n$, we have natural morphisms
$$\pi_n : {\cal L}(X) \rightarrow {\cal L}_n(X) \ \ {\rm and } \ \ \pi^m_n :
{\cal L}_m(X) \rightarrow {\cal L}_n(X),$$
obtained by truncation.  Note that ${\cal L}_0(X) = X$ and that ${\cal L}_1(X)$
is the tangent bundle of $X$.  For any arc $\gamma$ on $X$ (\textit{i.e.} a $K$-arc for
some field $K$ containing $k$), we call $\pi_0(\gamma)$ the origin of the
arc
$\gamma$.

\section{Motivic zeta function}
\subsection{}\label{not}
Let $k$ be a field and let $X$ be a smooth connected
variety over $k$.
Let $n \geq 1$ be an integer. 
The morphism $f : X \rightarrow \AA^1_k$
induces a morphism $f : {\cal L}_n(X) \rightarrow {\cal L}_n (\AA^1_k)$.

Any point $\alpha$ of ${\cal L}(\AA^1_k)$, resp. ${\cal L}_n(\AA^1_k)$, 
yields a $K$-rational point, for some field $K$ containing $k$, and
hence a power series $\alpha(t) \in K[[t]],$ resp. $\alpha(t) \in
K[[t]]/t^{n+1}$.
This yields maps
$${\rm ord}_t : {\cal L}(\AA^1_k) \rightarrow \NN \cup \{ \infty \},
\quad {\rm ord}_t : {\cal L}_n(\AA^1_k) \rightarrow \{ 0,1,\cdots,n,\infty
\},$$
with ${\rm ord}_t \alpha$ the largest $e$ such that $t^e$ divides $\alpha(t)$.

We set
$$\cX_n := \{ \varphi \in {\cal L}_n(X) \, | \, {\rm ord}_t f(\varphi) = n
\}.$$
This is a locally closed subvariety of ${\cal L}_n(X)$.  Note that $\cX_n$
is actually an $X_{0}$-variety, through the morphism $\pi^n_0 : {\cal L}_n(X)
\rightarrow X$.  Indeed $\pi^n_0(\cX_n) \subset X_0$, since $n \geq 1$. 
We consider the morphism
$$\bar f : \cX_n \rightarrow \GG_{m,k} := \AA^1_k \setminus \{
0 \},$$
sending a point $\varphi$ in $\cX_n$ to the coefficient of $t^n$ in
$f (\varphi)$.  There is a natural action of $\GG_{m,k}$ on $\cX_n$
given by $a \cdot \varphi(t) = \varphi(at)$, where $\varphi(t)$ is the vector
of
power series corresponding to $\varphi$ (in some local coordinate system).
Since
$\bar f (a \cdot \varphi) = a^n \bar f (\varphi)$ it follows
that $\bar f$ is a locally trivial fibration.

We denote by $\cX_{n}^{1}$ the fiber $\bar f^{-1}(1)$.  Note that the
action of $\GG_{m,k}$ on $\cX_n$ induces a good action of $\mu_n$ (and
hence of $\hat \mu$) on $\cX_{n}^{1}$.  Since $\bar f$ is a locally
trivial fibration, the $X_0$-variety $\cX_{n}^{1}$ and the action of $\mu_n$
on it, completely determines both the variety $\cX_n$ and the morphism
$$(\bar f, \pi^n_0) : \cX_n \rightarrow  \GG_{m,k} \times X_0.$$  Indeed it is easy to verify that $\cX_n$, as a
$(\GG_{m,k} \times X_0)$-variety, is isomorphic to $\cX_{n}^{1}
\times^{\mu_n} \GG_{m,k}$,
the quotient of 
$\cX_{n}^{1}
\times \GG_{m,k}$ under the $\mu_n$-action defined by $a(\varphi,b) = (a
\varphi,a^{-1}b)$.

The motivic zeta function of $f : X \rightarrow
\AA^1_k$, is the power series over $\mathcal{M}^{\hat \mu}_{X_0}$ defined by
$$Z_f(T) := \sum_{n \geq 1} \, [\cX_{n}^{1}/X_{0}, \hat \mu] \,
\LL^{-nr} \, T^n,$$ with $r$ the dimension of $X$.
When $k$ is of characteristic zero  
this a rational function of $T$ (cf. \cite{motivic}, \cite{barc}).

For $x$ a closed point of $X_0$, one defines similarly
$\cX_{n,x}$ and $\cX_{n,x}^{1}$,
by requiring
the arcs to have their origin in $x$.

The local motivic zeta function of $f$
at $x$ is the power series in $\mathcal{M}^{\hat \mu}_{x}$
defined by
$$Z_{f,x}(T) := \sum_{n \geq 1} \, [\cX_{n,x}^{1}, \hat \mu] \,
\LL^{-nr} \, T^n.$$

\begin{lem}\label{ho}Assume $X$ is the affine space $\AA^r_k$
and $f$ is a homogeneous polynomial of degree $d$ on $\AA^r_k$.
Then the equality
$$Z_f (T) =  \LL^r T^{-d} Z_{f,0} (T)$$ holds
in 
$\mathcal{M}^{\hat \mu}_{k} [[T]]$.
\end{lem}

\begin{proof}By homogeneity,
the mapping $\varphi \mapsto t \varphi$
induces an isomorphism between $\cX_n$ and $\pi_{n +1}^{n+d}(\cX_{n + d, 
0})$
compatible with the fibrations
$\cX_n \rightarrow \GG_{m,k}$ and $\cX_{n + r, 0} \rightarrow \GG_{m,k}$.
Hence $[\cX_{n}^{1}, \hat \mu]
=
[\cX_{n + d,0}^{1}, \hat \mu] \LL^{- (d - 1)r}$
and
the result follows.
\end{proof}

\section{Statement of the main result}
\subsection{}\label{tsc}Let $G$ be 
a connected reductive $k$-algebraic group 
and $\rho$ be a linear representation of $G$
on 
$\AA^m_k$.
Write
$m = r_1 + r_2$.
Then set $X_1 = M_{m, r_1}$ 
and let 
$G \times \SL_{r_1}$
act on $X_1$
by 
$$\rho_1 (g, g_1) x_1 = \rho (g) x_1 {}^{t}g_1.$$
Similarly
set $X_2 = M_{m, r_2}$
and let 
$G \times \SL_{r_2}$
act on $X_2$
by 
$$\rho_2 (g, g_2) x_2 = {}^t \rho (g)^{-1} x_2 {}^{t}g_2.$$
Assume $G_j = {\rm Im} \rho_j$
is an irreducible
regular prehomogeneous vector space for $j = 1, 2$ with
singular locus the irreducible hypersurface
$f_j = 0$.
One considers the quotient spaces
$X_i /\SL_{r_i}$ as embedded by
the Plücker embedding
$$
X_i /\SL_{r_i} \hookrightarrow V_i := \AA^{\binom{m}{r_i}}.
$$
There is a natural isomorphism
$V_1 \simeq V_2$ under which
the
quotient spaces $X_1 /\SL_{r_1}$ and
$X_2 /\SL_{r_2}$ may be naturally 
identified.
Hence we write $V$ for both $V_1$ and $V_2$, and 
up to multiplying one of
them by a nonzero constant factor,
one may assume 
that both $f_1$ and $f_2$
are the pullback of the same homogeneous
polynomial $f$ of degree $d$
in $V$, cf. \cite{class}. In particular,
${\deg} f_j = r_j d$, for $1 \leq j \leq 2$.

\begin{theorem}\label{MT}Let $(G_1, X_1)$
and
$(G_2, X_2)$ be irreducible regular
prehomogeneous vector spaces 
which are
related by a castling transformation. Then the  relations
\begin{equation}\label{lll}
Z_{f_1} (T)[{\rm SL}_{r_{2}, k}]
\prod_{1 \leq j \leq r_1} (1 -T^d \LL^{-j})
=
Z_{f_2} (T)
[{\rm SL}_{r_{1}, k}]\prod_{1 \leq j \leq r_2} (1 -T^d \LL^{-j})
\end{equation}
and 
\begin{equation}\label{lllbis}\begin{split}
Z_{f_1, 0} (T)\prod_{1 \leq j \leq r_1} (T^{-d} - \LL^{-j})&
\prod_{1 \leq j \leq r_2} (1 - \LL^{-j})
=\\
&Z_{f_2, 0} (T)
\prod_{1 \leq j \leq r_2} (T^{-d} - \LL^{-j})
\prod_{1 \leq j \leq r_1} (1 - \LL^{-j})\\
\end{split}
\end{equation}
hold
in 
$\mathcal{M}^{\hat \mu}_{k} [[T]]$, with
$[{\rm SL}_{r, k}]
=
\LL^{r^2 - 1} \prod_{2 \leq i \leq r} (1 - \LL^{-i})$.
\end{theorem}

\section{Proof of Theorem \ref{MT}}
\subsection{}Fix $1 \leq r \leq m$ and set
$Z = M_{m, r}$.
The $k$-group scheme ${\bf GL}_m (k [[t]])$ acts naturally
on 
$\cL (Z)$ on the left, while the
$k$-group scheme ${\bf SL}_{r} (k [[t]])$ acts naturally
on 
$\cL (Z)$ on the right. Here the group of $K$-points
of 
${\bf GL}_m (k [[t]])$, resp. ${\bf SL}_{r} (k [[t]])$,
is
${\GL_m} (K [[t]])$, resp. ${\SL_{r}} (K [[t]])$,
for $K$ a field containing $k$.
We shall
consider the quotient space
$h : Z \rightarrow Y := Z /\SL_{r}$ as embedded by
the Plücker embedding
$$
Z /\SL_{r} \hookrightarrow V := \AA^{\binom{m}{r}}.$$
We set $\cZ := \cL (Z)$ and $\cY := \cL (Y)$.
We also denote by 
$\cZ'$ the subset of $\cZ$ consisting of matrices
with a non zero minor of order $r$.
\subsection{}For every sequence
$\underline e := (e_1, \ldots, e_r)$ in $\NN^r$,
we denote by $\cZ_{\underline e}$ the subset of
$\cZ$ consisting of matrices of the form
$M A$ with $M$ in ${\GL_m} (K [[t]])$, 
and $A$ an upper triangular matrix 
with 
$t^{e_1}, \dots, t^{e_r}$ on the diagonal
and coefficients in $K [[t]]$, for some field $K$ containing $k$.
The set $\cZ'$ is the disjoint union of the
subsets
$\cZ_{\underline e}$, for
$\underline e $ in $\NN^r$.

\subsection{}For
$w$ in $W := \cS_r$, the permutation group on
$r$ letters, we consider
the set $\cV_{\underline e, w}$ of matrices
$x = (x_{i, j})_{1 \leq i \leq m, 1 \leq j \leq r}$ in $\cZ'$
of the form
$$
x_{i, j}
=
a_{i, j} u_j t^{e_j}
+
\sum_{1 \leq k < j} \lambda_{k, j} a_{i, k},
$$
with 
$\lambda_{k, j}$
and $a_{i, k}$ in $K [[t]]$,
$u_j$ a unit in $K [[t]]$,
and
$v (a_{i,j}) \geq 1$, for $i < w (j)$,
$a_{w (j), j} = 1$,
$a_{w (k), j} = 0$ for $k <j$.
We denote by 
$\cV_{\underline e}$
the disjoint union of the subsets $\cV_{\underline e, w}$
and
we denote by 
$\cW_{\underline e}$ the image of
$\cV_{\underline e}$ in $\cY$.
Let us remark that
$\cZ_{\underline e}$, resp. $h (\cZ_{\underline e})$,
is the union of a finite numbers
of 
${\bf GL}_m (k [[t]])$-translates of $\cV_{\underline e}$,
resp. $\cW_{\underline e}$.

\subsection{}Let us  recall the notion of piecewise trivial fibration.
Let $X$, $Y$ and $F$ be algebraic varieties over $k$,
and let
$A$, resp. $B$, be a constructible subset of $X$,
resp. $Y$. 
We say that
a map
$\pi : A \rightarrow B$ is a
piecewise trivial fibration with fiber
$F$, if there exists a finite partition of $B$ in subsets $S$ which are
locally closed
in $Y$ such that $\pi^{- 1} (S)$ is locally closed in $X$ and
isomorphic, as
a variety over $k$, to $S \times F$, with $\pi$
corresponding under the isomorphism to the projection
$S \times F \rightarrow S$.

\begin{lem}\label{ML}
For $n \geq |\underline e|$,
with
$|\underline e| :=
\sum_{1 \leq i \leq r} e_i$,
the canonical morphism
$$
\pi_n(\cV_{\underline e, w})
\longrightarrow
\pi_n (\cW_{\underline e})
$$
is a piecewise trivial fibration with fiber
$$
Z_w \times 
\AA^{n (r^2- 1) + |\underline e| ((m- r)r +1)}_k
\prod_{1 \leq i \leq r} \AA^{- (m +1 -i)e_i}_k,
$$
with
$Z_w := (\GG_{m, k})^{r -1} \prod_{1 \leq i \leq r} \AA^{r - 1- m_i}_k$
and $m_i$ the number of integers $1 \leq k < w (i)$
with $k \not= w (j)$ for $j <i$.
\end{lem}

\begin{proof}Let us first consider the case where $w ={\rm id}$.
Let
$x$ be an $(m, r)$
matrix in $\cV_{\underline e, {\rm id}}$.
For $r + 1 \leq k \leq m$ and $j \leq r$, 
we
shall denote by
$\Delta^{k, j}$ 
the determinant of the $(r, r)$ submatrix
of $x$ obtained by removing the $j$-th line from the
$(r +1, r)$ submatrix of $x$ obtained by keeping the first $r$ lines
together with the $k$-th line.
We can take $\Delta^{k, j}$, for $r + 1 \leq k \leq m$ and $1 \leq j \leq r$, 
as coordinates on 
$\cW_{\underline e}$.
Set $\Delta := \prod_{1 \leq i \leq r} u_i \pi^{e_i}$.
We have the relations
$$
\Delta^{k, r} = a_{k, r} \Delta,
$$
$$
\Delta^{k, r- 1} = a_{r, r- 1} \Delta^{k, r} - a_{k, r - 1} \Delta,
$$
$$ \dots$$
$$
\Delta^{k, i} = \sum_{1 \leq j \leq r - i} (-1)^{j + 1} a_{i +j, i}
\Delta^{k, i+j}
+ (-1)^{r -i} a_{k, i} \Delta,
$$
$$ \dots$$
$$
\Delta^{k, 1} = a_{2, 1} \Delta^{k, 2} - a_{3,1} \Delta^{k, 3}
+ \ldots + (-1)^r a_{r, 1} \Delta^{k, r} +
(-1)^{r + 1} a_{k, 1} \Delta,
$$
from which the result follows.
For general $w$, the situation reduces to the former one,
since, up to renumbering the rows, the situation is just the same
except that some $a_{i, j}$'s are now required to be of valuation
$\geq 1$, which has the sole effect of replacing $Z_{\rm id}$ by
$Z_w$
in the statement.
\end{proof}

\begin{lem}\label{ut}Set $\cY_{0} := h (\cZ_{\underline 0})$
and, for $k \geq 0$, 
$\cY_{k} = t^k \cY_{0}$.
\begin{enumerate}
\item[(1)]For every $\underline e$ in $\NN^r$,
$h (\cZ_{\underline e})  = \cY_{\vert \underline e \vert}$.
\item[(2)]Assume $n \geq k \geq 0$.
For every
constructible subset 
$A$ of
$\pi_{n} (\cY_{k})$, 
consider the preimage $[h^{-1} (A)]$ of $A$ in $\pi_{n} (\cZ)$.
The relation
$$
[h^{-1} (A)]
=
[A]
[\SL_{r, k}]
\LL^{n (r^2 -1) + k ((m - r)r + 1)}
\sum_{\vert \underline e \vert = k}
\prod_{1 \leq i \leq r} \LL^{- (m + 1- i)e_{i}}
$$
holds in $\cM_k$. Furthermore,
if $A$ is constructible with
$\hat \mu$-action,
the same relation still
holds in 
$\mathcal{M}^{\hat \mu}_{k}$.
\end{enumerate}
\end{lem}

\begin{proof}Follows directly from
Lemma \ref{ML} and the relation
$$\sum_{w \in W} [Z_w]
=
\LL^{r^2 - 1} \prod_{2 \leq i \leq r} (1 - \LL^{-i})
=
[{\rm SL}_{r, k}].$$
\end{proof}

\subsection{}Let $f$ be a homogeneous polynomial of degree $d$
on $V$.
We set
$$\cY_{n, k} :=
\{ \varphi \in {\cal L}_n(Y) \, | \, 
\varphi \in \pi_{n} (\cY_{k}) \quad \text{and} \quad
{\rm ord}_t f (\varphi) = n
\}.$$
As in \ref{not}, we define $\bar f : \cY_{n, k} \rightarrow
\GG_{m, k}$, and we
set
$\cY_{n, k}^{1} := \bar f^{-1} (1)$, which is a variety with
$\hat \mu$-action.
We shall also consider the varieties
$\cX_{n, k}$
and $\cX_{n, k}^1$ which are the inverse images of respectively
$\cY_{n, k}$
and $\cY_{n, k}^1$ in $\cL_{n} (Z)$.

\begin{lem}\label{hom}The relation
$$
[\cY_{n, k}^1]
=
[\cY_{n- kd, 0}^1] \LL^{k (d - 1) \dim Y}
$$ holds
in $\mathcal{M}^{\hat \mu}_{k}$.
\end{lem}

\begin{proof}By homogeneity,
the mapping $y \mapsto t^k y$
induces an isomorphism between $\cY_{n - kd, 0}$ and
$\pi_{n - k(d - 1)}^{n}(\cY_{n, k})$
compatible with the fibrations
onto $\GG_{m,k}$.
\end{proof}

\subsection{}Let us now consider
the motivic zeta function of $f \circ h$
(or more precisely its image in
$\mathcal{M}^{\hat \mu}_{k} [[T]]$),
$$Z_{f \circ h}
(T) := \sum_{n \geq 1} \, [\cX_{n}^{1}, \hat \mu] \,
\LL^{-n \dim Z} \, T^{n},$$
where $\cX_{n}^{1}$
is the disjoint union of
the varieties with $\hat \mu$-action
$\cX_{n, k}^{1}$, $k$ in $\NN$.
Remark that for $kd > n$, $\cX_{n, k}$ and
$\cY_{n, k}$ are empty, hence for fixed $n$
only a finite number of the $\cX_{n, k}^{1}$'s are non empty.
We shall also consider the power series
$$
Z_{f \circ h}^{0} (T)
:= \sum_{n \geq 1} \, [\cX_{n, 0}^{1}, \hat \mu] \,
\LL^{-n \dim Z} \, T^n
$$
and
$$
Z_{f}^0 (T) :=
\sum_{n \geq 1} \, [\cY_{n, 0}^{1}, \hat \mu] \,
\LL^{-n \dim Y} \, T^n
$$
in
$\mathcal{M}^{\hat \mu}_{k} [[T]]$.

\begin{prop}The  relations 
\begin{equation}\label{e1}
Z_{f \circ h} (T)
=
\prod_{1 \leq i \leq r}
(1- \LL^{- (m + 1 -i)} T^d)^{-1}
Z_{f \circ h}^{0} (T)
\end{equation}
and
\begin{equation}\label{e2}
Z_{f \circ h} (T)
= [\SL_{r, k}]
\prod_{1 \leq i \leq r}
(1- \LL^{- (m + 1 -i)} T^d)^{-1}
Z_{f}^{0} (T)
\end{equation}
hold
in
$\mathcal{M}^{\hat \mu}_{k} [[T]]$.
\end{prop}

\begin{proof}Since, as follows from Lemma \ref{ut},
$$
Z_{f \circ h}^{0} (T)
= [\SL_{r, k}]
Z_{f}^{0} (T),
$$
it is enough to prove (\ref{e2}).
By Lemma  \ref{ut}, we
have
$$
[\cX_{n, k}^{1}]
\LL^{-n \dim Z}
=
[\SL_{r, k}]
[\cY_{n, k}^{1}]
\LL^{-n \dim Y}
\LL^{k \dim Y}
\sum_{\vert \underline e \vert = k}
\prod_{1 \leq i \leq r} \LL^{-  (m + 1 - i) e_{i}}.
$$
It follows from
Lemma \ref{hom} that we may rewrite the last equality as
$$
[\cX_{n, k}^{1}]
\LL^{-n \dim Z}
=
[\SL_{r, k}]
[\cY_{n -kd, 0}^{1}]
\LL^{- (n -kd) \dim Y}
\sum_{\vert \underline e \vert = k}
\prod_{1 \leq i \leq r} \LL^{-  (m + 1 - i) e_{i}},
$$
and we get the result by summing up the series.
\end{proof}

\subsection{}Now we can conclude the
proof of Theorem \ref{MT}. Relation (\ref{lll})
follows from writing
(\ref{e2}) for both $X_{1}$ and $X_{2}$, 
and 
(\ref{lllbis}) follows from (\ref{lll}) together with Lemma \ref{ho}.
\qed

\section{Applications to the Milnor fibre and the Hodge spectrum}\label{jj}

\subsection{Monodromy}\label{mono}In this subsection \ref{mono}
we assume that $k = \CC$.
Let $x$ be a point of
$X_0 = f^{-1}(0)$.  We fix a smooth metric on $X$.
We set $X^\times_{\varepsilon,\eta} := B(x,\varepsilon) \cap
f^{-1}(D^\times_\eta)$, with $B(x,\varepsilon)$ the open ball of radius $\varepsilon$
centered at $x$ and $D^\times_\eta := D_\eta \setminus \{ 0 \}$, with $D_\eta$
the open disk of radius $\eta$ centered at $0$.  For $0 < \eta \ll \varepsilon \ll
1$,
the restriction of $f$ to $X^\times_{\varepsilon,\eta}$ is a locally trivial
fibration, called the Milnor fibration, onto $D^\times_\eta$ with fiber
$F_x$, the Milnor fiber at $x$.  The action of a characteristic
homeomorphism of this fibration on cohomology gives rise to the monodromy
operator
$$M_x : H^\cdot (F_x,\QQ) \rightarrow H^\cdot (F_x,\QQ).$$

\subsection{The motivic Milnor fibre}\label{mmf}
It is a remarkable fact,
proved in
\cite{motivic} and
\cite{barc}, that
expanding the rational function $Z_f(T)$ as a power series in $T^{-1}$ and
taking
minus its constant term, yields a well defined element 
$\cS_f = 
- \lim_{T \rightarrow  \infty}Z_f(T)$
of $\mathcal
{M}_{X_0}^{\hat
\mu}$.
This follows from 
a formula for the motivic zeta function in terms of embedded resolutions of
$X_0$ which we now recall.

Let $(Y,h)$ be a
resolution of $f$.  By this, we mean that $Y$ is a
nonsingular and irreducible algebraic variety over $k, h : Y \rightarrow X$ is
a
proper morphism, that the restriction $h : Y \setminus h^{-1}(X_0) \rightarrow
X
\setminus X_0$ is an isomorphism, and that $h^{-1}(X_0)$ has only normal
crossings
as a subvariety of $Y$. 

We  denote by $E_i, i \in J$, the irreducible  components (over $k$) of
$h^{-1}(X_0)$.  For each $i \in J$, denote by $N_i$ the multiplicity of $E_i$
in
the divisor of $f \circ h$ on $Y$, and by $\nu_i - 1$ the multiplicity of $E_i$
in the divisor of $h^\ast dx$, where $dx$ is a local non vanishing volume form
at any point of $h(E_i)$, i.e. a local generator of the sheaf of differential
forms of maximal
degree.  For $i \in J$ and $I \subset J$, we consider the nonsingular varieties
$E^\circ_i := E_i \setminus \cup_{j \ne i} E_j, E_I = \cap_{i \in I} E_i$, and
$E^\circ_I := E_I \setminus \cup_{j \in J \setminus I} E_j$.

Let $m_I = {\gcd}(N_i)_{i \in I}$.  We introduce an unramified
Galois cover $\tilde E^\circ_I$ of $E^\circ_I$, with Galois group $\mu_{m_I}$,
as follows.  Let $U$ be an affine Zariski open subset of $Y$, such that, on
$U$, $f \circ h = uv^{m_I}$, with $u$ a unit on $U$ and $v$ a morphism from $U$
to $\AA^1_k$.  Then the restriction of $\tilde E^\circ_I$ above $E^\circ_I
\cap U$, denoted by $\tilde E_I^\circ \cap U$, is defined as
$$\{ (z,y) \in \AA^1_k \times (E^\circ_I \cap U) | z^{m_I} = u^{-1} \}.$$
Note that $E^\circ_I$ can be covered by such affine open subsets $U$ of $Y$.
Gluing together the covers $\tilde E^\circ_I \cap U$, in the obvious way, we
obtain the cover $\tilde E^\circ_I$ of $E^\circ_I$ which has a natural
$\mu_{m_I}$-action (obtained by multiplying the $z$-coordinate with the elements
of $\mu_{m_I}$).  This $\mu_{m_I}$-action on $\tilde E^\circ_I$ induces an
$\hat \mu$-action on $\tilde E^\circ_I$ in the obvious way.

\begin{theorem}[\cite{lef}, \cite{looi}]\label{3.3.1}
With the previous notations, the
following relation
holds in $\mathcal{M}^{\hat \mu}_{X_0} [[T]]$ :
$$Z_f(T) =  \sum_{\emptyset \ne I \subset J} \, (\LL - 1)^{|I|-1} \,
[ \tilde
E^\circ_I/X_0, \hat \mu] \,
\prod_{i
\in I} \frac{\LL^{-\nu_i} T^{N_i}}{1 - \LL^{-\nu_i} T^{N_i}}.$$
\end{theorem}

\begin{definition}[\cite{motivic}, \cite{lef}, \cite{barc}]\label{3.5.3}
Expanding the rational function $Z_f(T)$ as a power series in $T^{-1}$
and taking
minus its constant term, yields a well defined element of $\mathcal
{M}_{X_0}^{\hat
\mu}$, namely
$$\mathcal{S}_f := - \lim_{T \rightarrow \infty} Z_f(T) :=
\sum_{\emptyset \ne I
\subset J} (1 - \LL)^{|I|-1} [ \tilde E^\circ_I].$$
Moreover we set $\mathcal{S}_{f, x} := {\rm Fiber}_x(\mathcal{S}_f)$ in
$\mathcal{M}^{\hat
\mu}_k$.
We also have $\cS_{f, x} = 
- \lim_{T \rightarrow  \infty}Z_{f, x}(T)$.
\end{definition}

There is strong evidence
that $\mathcal{S}_{f, x}$ is the correct virtual motivic
incarnation of
the Milnor fiber $F_x$ of $f$ at $x$ (which is in itself not at all motivic).

In the present setting we deduce from Theorem \ref{MT},
the following relation between the virtual motivic Milnor fibres
at the origin of two irreducible regular
prehomogeneous vector spaces 
which are
related by a castling transformation.

\begin{theorem}\label{milnor}Assume $k$ is a field of characteristic zero.
Let $(G_1, X_1)$
and
$(G_2, X_2)$ be irreducible regular
prehomogeneous vector spaces 
which are
related by a castling transformation. Then the  relation
\begin{equation}\label{mmm}
\cS_{f_1, 0} (T)
\prod_{1 \leq j \leq r_2} (1 - \LL^{j})
=
\cS_{f_2, 0} (T)
\prod_{1 \leq j \leq r_1} (1 - \LL^{j})
\end{equation}
holds
in 
$\mathcal{M}^{\hat \mu}_{k}$.
\end{theorem}

\begin{proof}Follows directly from (\ref{lllbis}).
\end{proof}

\subsection{Hodge structures}\label{3.1.2}From now on in the remaining
of this
section \ref{jj}, we shall assume $k = \CC$.
A Hodge structure is a finite dimensional
$\QQ$-vector
spa\-ce $H$ together with a bigrading $H \otimes \CC =
\oplus_{p,q \in \ZZ} H^{p,q}$, such that $H^{q,p}$ is the complex conjugate
of $H^{p,q}$  and each weight $m$ summand, $\oplus_{p+q=m} H^{p,q}$, is
defined
over $\QQ$.  The Hodge structures, with the evident notion of morphism, form
an abelian category $\HS$ with tensor product.  
The elements of the Grothendieck
group $K_0 (\HS)$ of this abelian category are representable as a formal
difference of Hodge structures $[H] - [H^\prime]$, and  $[H] = [H^\prime]$ iff 
$H \cong H^\prime$.  Note that $K_0 (\HS)$ becomes a ring with respect to the
tensor product.

A mixed Hodge structure is a finite dimensional $\QQ$-vector space $V$
with
a finite increasing filtration $W_\bullet V$, called the weight
filtration, such
that the associated graded vector space ${\rm Gr}^W_\bullet(V)$ underlies a Hodge
structure having ${\rm Gr}^W_m(V)$ as weight $m$ summand.  Note that $V$ determines
in a natural way an element $[V]$ in $K_0 (\HS)$, namely $[V] := \sum_m
[{\rm Gr}^W_m(V)]$.

When $X$ is an algebraic variety over $k = \CC$, the simplicial cohomology
groups $H^i_c(X, \QQ)$ of $X$, with compact support, underly a natural mixed
Hodge structure, and the Hodge characteristic $\chi_h(X)$ of $X$ (with
compact support) is defined by
$$\chi_h(X) := \sum_i (-1)^i [H^i_c(X,\QQ)] \in K_0(\HS).$$
This yields a map $\chi_h :
\Var_{\CC} \rightarrow K_0 (\HS)$
which factors through $\mathcal{M}_k$,
because $\chi_h (\AA^1_{k})$ is actually invertible in the ring
$K_0 (\HS)$.  When
$X$ is proper and smooth, the mixed Hodge structure on $H^i_c(X,\QQ)$ is in
fact a Hodge structure, the weight filtration being concentrated in weight $i$.

\subsection{The Hodge spectrum}\label{3.1.3}
The cohomology groups $H^i (F_x,\QQ)$ of
the Milnor fiber $F_x$ carry a natural mixed Hodge structure 
(\cite{St1}, \cite{Va}, \cite{Sa1}, \cite{Sa2}),
which
is compatible with the semi-simplification
of the monodromy operator $M_x$.  Hence we can define the Hodge
characteristic $\chi_h(F_x)$ of $F_x$ by
$$\chi_h(F_x) := \sum_i (-1)^i [H^i(F_x,\QQ)] \in K_0(\HS).$$
Actually by taking into account the monodromy action we can consider
$\chi_h(F_x)$ as an element of the Grothendieck group $K_0(\HS^{\rm mon})$ of
the abelian category $\HS^{\rm mon}$ of Hodge structures 
with an endomorphism of finite order.
Again $K_0(\HS^{\rm mon})$ is a ring by the tensor product.

There is a natural linear map, called the Hodge spectrum
$$\hsp :
K_0(\HS^{\rm mon }) \rightarrow \ZZ[t^{1/\ZZ}] := \cup_{n
\geq 1} \ZZ[t^{1/n},t^{-1/n}],$$
with $\hsp ([H]) := \sum_{\alpha \in \QQ \cap[0,1[} t^\alpha (\sum_{p,q
\in
\ZZ} \dim(H^{p,q})_\alpha)t^p$, for any Hodge structure $H$ with an
endomorphism of finite order, where $H^{p,q}_\alpha$ is the generalized
eigenspace of $H^{p,q}$ with respect to the eigenvalue
$e^{2\pi\sqrt{-1}\alpha}$.

We recall that $\hsp (f,x) := (-1)^{d - 1}\hsp (\chi_h(F_x)-1)$ is called the Hodge
spectrum of
$f$ at $x$.

We denote by $\chi_h$ the canonical ring
homomorphism (called the Hodge characteristic)
$$\chi_h : \mathcal{M}^{\hat \mu}_k \rightarrow K_0({\HS^{\rm mon}}),$$
which associates to any complex algebraic variety $Z$, with a good
$\mu_n$-action, its Hodge characteristic  together with the
endomorphism induced by $Z \rightarrow Z : z \mapsto e^{2 \pi\sqrt{-1}/n}z$.

\begin{theorem}[\cite{motivic}]\label{3.5.5}
Assume the above notation
with $k
= \CC$.  Then we have the following equality in
$K_0({\HS^{\rm mon}})$ :
$$\chi_h(F_x) = \chi_h(\mathcal{S}_{f, x}).$$
\end{theorem}

In particular it follows that
$$
\hsp (f,x) := (-1)^{d - 1}\hsp (\chi_h(\mathcal{S}_{f, x})-1),
$$
thus the motivic zeta function
$Z_{f, x}(T)$ completely determines the Hodge spectrum of $f$ at $x$.

Hence we deduce from Theorem \ref{milnor}, since
$\hsp (\chi_h (\LL)) = t$,
the following:

\begin{cor}\label{sp}Assume $k = \CC$.
Let $(G_1, X_1)$
and
$(G_2, X_2)$ be irreducible regular
prehomogeneous vector spaces 
which are
related by a castling transformation. Then the following relation
holds between the Hodge spectra of $f_1$ and $f_2$ at $0$
\begin{equation*}
\frac{1 + (-1)^{m r_1 - 1}\hsp (f_1, 0)}{\prod_{1 \leq j \leq r_1} (1 - t^j)}
=
\frac{1 + (-1)^{m r_2 - 1}\hsp (f_2, 0)}{\prod_{1 \leq j \leq r_2} (1 - t^j)}.
\qed
\end{equation*}
\end{cor}

\section{The general  case}\label{gen}

\subsection{}We now consider the general
case of prehomogeneous vector spaces
which are not necessarily  irreducible and regular.
In this setting we are given
a connected linear algebraic group
$G$ together
with a linear action
$\rho : G \rightarrow
{\rm End} X$
on a finite dimensional (linear) affine space
$X$ over $k$, such that the action
is transitive on
a dense open set $O = X \setminus S$.
Let $S_j$, $1 \leq j \leq \ell$, be the $k$-irreducible components of
the singular set
$S$ which are of codimension 1 in $X$
and choose for every $j$ a defining equation $f_j = 0$ of $S_j$.
Recall that a non-zero $k$-rational function on $X$ is
a $k$-relative invariant of the $G$-action on $X$,
if there exists a $k$-rational character $\nu$
of $G$ such that
$$
f (\rho (g) x) = \nu (g) f (x) 
$$
for every $g$ in $G$ and $x$ in $X$.
The functions $f_j$ are 
$k$-relative invariants and furthermore 
they are a basis of $k$-relative invariants in the sense that
any
$k$-relative invariant of $(G, \rho)$ is of the form
$c f_1^{\mu_1} \dots  f_{\ell}^{\mu_{\ell}}$
with $c$ in $k^{\times}$ and $\mu_j$ in $\ZZ$.

\subsection{Motivic zeta function for several functions}\label{jkj}
Let $X$ be a smooth $k$-variety of dimension $d$
and consider $\ell$ functions $f_j : X \rightarrow \AA^1_k$.
We set $X_0 := \cap_{1 \leq j \leq \ell} (f_j = 0)$.

For every $n = (n_1, \dots, n_{\ell})$
in $(\NN^{\times})^{\ell}$, we
set $|n| := \sum_{1 \leq j \leq \ell} n_j$ and
define
$$\cX_n := \Bigl\{x \in \cL_{|n|} (X) \Bigm|
\ord_t f_j = n_j, 1 \leq j \leq \ell\Bigr\}.$$
Similarly as in
\ref{not}, we have a natural morphism
$\bar f : \cX_n \rightarrow \GG_{m, k}^{\ell}$,
which makes $\cX_n$ a $X_0 \times \GG_{m, k}^{\ell}$-variety,
for
$n$
in $(\NN^{\times})^{\ell}$.
The motivic 
zeta function attached to $f = (f_1, \dots, f_{\ell})$
is the formal series
$$
Z_f (T) :=
\sum_{n \in (\NN^{\times})^{\ell}}
[\cX_n /
X_0 \times \GG_{m, k}^{\ell}] \LL^{- |n| d} T^n
$$
in $\cM_{X_0 \times \GG_{m, k}^{\ell}}[[T_1, \dots T_{\ell}]]$.
When $k$ is of characteristic zero $Z_f$ is a rational function of
$T = (T_1, \dots T_{\ell})$ (cf. \cite{motivic}, \cite{looi}, \cite{barc}).
When $\ell = 1$, the relation with the definition in 
\ref{not} is the following:
the 
zeta function we just defined is the image of the former
one by the canonical morphism
$\cM^{\hat \mu}_{X_0} \rightarrow
\cM_{X_0 \times \GG_{m, k}}$ which to the class
of a $X_0$-variety $Y$
with $\hat \mu$-action assigns the class of
$Y \times^{\hat \mu} \GG_{m, k}$ in 
$\cM_{X_0 \times \GG_{m, k}}$.
For $x$ a point in  $X_0 (k)$,
one defines similarly as in the $\ell = 1$ case,
$Z_{f, x} (T)$
in $\cM_{\GG_{m, k}^{\ell}}[[T_1, \dots T_{\ell}]]$.

\subsection{}
As in 
\ref{tsc}
we consider  a connected linear algebraic group $G$
over $k$ and $\rho$
be a representation of $G$ in
$\AA^m_k$.
Write
$m = r_1 + r_2$
and define
$X_j$, $G_j$ and $\rho_j$, $j = 1, 2$ as in \ref{tsc}.
We assume they give rise to
prehomogeneous vector spaces for $j = 1, 2$,
and we say that
the prehomogeneous vector spaces $(G_1, X_1)$
and
$(G_2, X_2)$ are related by a castling transformation.
We choose a basis of $k$-relative invariants
$f_1 = (f_{1, i})$, $1 \leq j \leq \ell_1$ and
$f_2 = (f_{1, 2})$, $1 \leq j \leq \ell_2$, respectively.
As in \ref{tsc},
we
consider the quotient spaces
$X_i /\SL_{r_i}$ as embedded by
the Plücker embedding
the same affine space $V$, and up to renumbering
the
$f_{1, i}$'s and multiplying them by non-zero
constants,
one may assume 
that both  the
$f_{1, i}$'s and the
$f_{2, i}$'s
are the pullback of the same homogeneous
polynomials $f_i$ of degree $d_i$
in $V$, cf. \cite{class}, \cite{parabolic}.
In particular, we may write $\ell_1 = \ell_2 = \ell$.
We set $d = (d_i)$ in $(\NN^{\times})^{\ell}$.

We have the following generalization of Theorem \ref{MT}.

\begin{theorem}\label{MTbis}Let $(G_1, X_1)$
and
$(G_2, X_2)$ be 
prehomogeneous vector spaces 
related by a castling transformation. 
Then the  relations
\begin{equation}\label{www}
Z_{f_1} (T)[{\rm SL}_{r_{2}, k}]
\prod_{1 \leq j \leq r_1} (1 -T^d \LL^{-j})
=
Z_{f_2} (T)
[{\rm SL}_{r_{1}, k}]\prod_{1 \leq j \leq r_2} (1 -T^d \LL^{-j})
\end{equation}
and 
\begin{equation}\label{wwwbis}
\begin{split}
Z_{f_1, 0} (T)\prod_{1 \leq j \leq r_1} (T^{-d} - \LL^{-j})
&\prod_{1 \leq j \leq r_2} (1 - \LL^{-j})
=\\
&Z_{f_2, 0} (T)
\prod_{1 \leq j \leq r_2} (T^{-d} - \LL^{-j})
\prod_{1 \leq j \leq r_1} (1 - \LL^{-j})\\
\end{split}
\end{equation}
hold
in 
$\mathcal{M}_{\GG_{m, k}^{\ell}} [[T]]$, 
with $T = (T_1, \dots,
T_{\ell})$ and 
$[{\rm SL}_{r, k}]
=
\LL^{r^2 - 1} \prod_{2 \leq i \leq r} (1 - \LL^{-i})$.
\end{theorem}

\begin{proof}The proof is just the same as the proof of Theorem \ref{MT}.
\end{proof}

\begin{remark}There is an obvious generalization of Theorem \ref{MTbis}
to parabolic castling transforms as
introduced in \cite{ter}, cf. \cite{parabolic}. Details are left to the reader.
\end{remark}

\subsection{}Similarly as for the case of one function 
in \ref{3.5.3},
Guibert  proves in \cite{guibert} that, for every $\alpha$ in
$(\NN^{\times})^{\ell}$,
$- \lim_{T \rightarrow  \infty}Z_f(T^{\alpha})$
is a well defined
element of $\cM_{X_0 \times \GG_{m, \CC}^{\ell}}$, independent of
$\alpha$. Let us denote it by $\cS_f$
and set $\cS_{f, x} := {\rm Fiber}_x (\cS_f)$. We also have
$\cS_{f, x} = - \lim_{T \rightarrow  \infty}Z_{f, x}(T^{\alpha})$,
for every $\alpha$ in
$(\NN^{\times})^{\ell}$.

\begin{remark}When $k = \CC$,
G. Guibert shows in 
\cite{guibert}
how Sabbah's Alexander invariants
of $f$ (cf. \cite{sabbah})
may be recovered from the motivic zeta function $Z_f$.
\end{remark}
\begin{remark}It would be interesting to investigate
what information is contained in the Hodge realization
of 
$\cS_f$.
It is quite likely that there should be some connections
with
recent work A. Libgober in
\cite{lib}.
\end{remark}

The following statement follows directly from
Theorem \ref{MTbis}.

\begin{theorem}\label{milnorbis}Let $(G_1, X_1)$
and
$(G_2, X_2)$ be 
prehomogeneous vector spaces 
related by a castling transformation. Then the  relation
\begin{equation}\label{mmmbis}
\cS_{f_1, 0} (T)
\prod_{1 \leq j \leq r_2} (1 - \LL^{j})
=
\cS_{f_2, 0} (T)
\prod_{1 \leq j \leq r_1} (1 - \LL^{j})
\end{equation}
holds
in 
$\cM_{X_0 \times \GG_{m, \CC}^{\ell}}$.
\end{theorem}

\bibliographystyle{amsplain}

\begin{thebibliography}{SGA}

\bibitem{motivic}
J. Denef, F. Loeser,
\textit{Motivic Igusa zeta functions}, 
J. Algebraic Geom. \textbf{7} (1998), 
505--537. 


\bibitem{lef}J. Denef, F. Loeser,
\textit{Lefschetz numbers of iterates of the monodromy and truncated arcs},
to appear in Topology, available at math.AG/0001105.

\bibitem{barc}
J. Denef, F. Loeser, \textit{Geometry on arc spaces of algebraic varieties},
to appear in
Proceedings of the Third European Congress of Mathematics, Barcelona 2000
available at math.AG/0006050.


\bibitem{guibert}
G. Guibert,
\textit{Fonction zêta motivique associée à une famille de séries
de deux variables},
C. R. Acad. Sci. Paris Sér. I Math. 
\textbf{333} (2001), 457--460. 





\bibitem{parabolic}
H. Hosokawa, \textit{Igusa local zeta functions
and parabolic castling transformation of
prehomogeneous vector spaces}, J. Number Theory \textbf{74} (1999), 
148--171.

\bibitem{igusa0}
J. Igusa, 
\textit{On functional equations of complex powers},
Invent. Math. \textbf{85} (1986), 1--29. 





\bibitem{igusa}
J. Igusa, 
\textit{On the arithmetic of a singular invariant},
Amer. J. Math. \textbf{110} (1988), 197--233.

\bibitem{aa}
J. Igusa,
\textit{$b$-functions and $p$-adic integrals}, Algebraic analysis, Vol. I, 231--241,
Academic Press, Boston, MA, 1988. 


\bibitem{igusaclass}
J. Igusa, 
\textit{Local zeta functions of certain prehomogeneous vector spaces},
Amer. J. Math. \textbf{114} (1992), 251--296. 

\bibitem{kimura}
T. Kimura, 
\textit{The $b$-functions and holonomy diagrams of irreducible regular prehomogeneous vector spaces},
Nagoya Math. J. \textbf{85} (1982), 1--80. 
\bibitem{lib}
A. Libgober,
\textit{Hodge decomposition of Alexander invariants},
preprint, October 2001, available at math.AG/0108018.

\bibitem{looi}
E. Looijenga, \textit{Motivic Measures}, in S{\'e}minaire Bourbaki,
expos{\'e} 874, Mars 2000, available at math.AG/0006220.

\bibitem{sabbah}
C. Sabbah,
\textit{Modules d'Alexander et $\cD$-modules},
Duke Math. J. \textbf{60} (1990), 729--814.


\bibitem{Sa1}M. Saito, \textit{Modules de Hodge polarisables}, Publ. Res.
Inst. Math. Sci., \textbf{24} (1988), 849--995.


\bibitem{Sa2}M. Saito, \textit{Mixed Hodge modules}, Publ.
Res. Inst. Math. Sci.,
\textbf{26} (1990), 221--333.



\bibitem{class}
M. Sato, T. Kimura,
\textit{A classification of irreducible prehomogeneous vector spaces
and their relative
invariants}, Nagoya Math. J. \textbf{65} (1977), 1--155. 



\bibitem{sato}
M. Sato, T. Shintani, 
\textit{On zeta functions associated with prehomogeneous vector spaces},
Ann. of Math. \textbf{100} (1974), 131--170.

\bibitem{ter}
Y. Teranishi, 
\textit{Relative invariants and
$b$-functions of prehomogeneous vector spaces
$(G\times{\rm GL}(d\sb
1,\cdots,d\sb r),\tilde p\sb 1, M(n,C))$},
Nagoya Math. J. \textbf{98} (1985), 139--156. 




\bibitem{St1}J. Steenbrink,
\textit{Mixed Hodge structures on the vanishing cohomology},
in Real and Complex Singularities, Sijthoff and Noordhoff, Alphen aan 
den Rijn,
1977, 525--563.
\bibitem{St2}J. Steenbrink, 
\textit{The spectrum of hypersurface
singularities}, in Th{\'e}orie de Hodge, Luminy 1987, Ast{\'e}risque, {\bf 179-180}
(1989), 163--184.



\bibitem{Va}A. Varchenko, \textit{Asymptotic Hodge structure in the vanishing
cohomology}, Math. USSR Izvestija, \textbf{18} (1982), 469--512.



\end{thebibliography}

\end{document}